\pdfoutput=1
\documentclass[journal]{IEEEtran}

\usepackage{amsmath,amssymb,mathrsfs,bm}
\usepackage{booktabs}
\usepackage{url}
\usepackage{cite}
\usepackage{graphicx}
\usepackage{xcolor}

\newcommand{\R}{\mathbb{R}}
\newcommand{\C}{\mathbb{C}}
\newcommand{\I}{\mathsf I}
\newcommand{\Pn}{\mathsf P_\nu}
\newcommand{\Pt}{\mathsf P_t}
\newcommand{\Sop}{\mathsf S}
\newcommand{\Kop}{\mathsf K}
\newcommand{\Kp}{\mathsf K'}
\newcommand{\Top}{\mathsf T}
\newcommand{\Aop}{\mathsf A}
\newcommand{\Bop}{\mathsf B}
\newcommand{\de}{\,\mathrm d}
\newcommand{\e}{\mathrm e}
\newcommand{\im}{\mathrm i}

\begin{document}

\title{Electric Combined-Field-Only P1 Galerkin BEM\\for PEC Scattering with Scalar Helmholtz Operators}

\author{Carlos P\'erez-Arancibia and Catalin Turc%
\thanks{C. P\'erez-Arancibia is with the Department of Applied Mathematics, University of Twente, Enschede, The Netherlands (e-mail: \protect\url{c.a.perezarancibia@utwente.nl}).}%
\thanks{C. Turc is with the Department of Mathematical Sciences, New Jersey Institute of Technology, Newark, NJ, USA (e-mail: \protect\url{catalin.turc@njit.edu}).}}

\markboth{P\'erez-Arancibia and Turc: Electric Combined-Field-Only P1 Galerkin BEM}{P\'erez-Arancibia and Turc: Electric Combined-Field-Only P1 Galerkin BEM}

\maketitle

\begin{abstract}
We develop a Galerkin boundary element method for the indirect electric combined-field-only integral equation (ECFOIE) for perfectly electrically conducting obstacles  with smooth surfaces. The method is free of low-frequency breakdown and of spurious resonances. It applies scalar Helmholtz operators componentwise, so that fast scalar solvers can be reused, with less memory than a Rao--Wilton--Glisson (RWG) discretization. The iteration count of the Calder\'on-regularized ECFOIE (R-ECFOIE) is constant under mesh refinement and, with a complex-wavenumber Calder\'on regularizer, grows only mildly with frequency. The discretization uses continuous, nodal piecewise-linear basis functions on flat triangular meshes, with essentially as many unknowns as an RWG discretization, and tangential and normal projections and the mean curvature of the obstacle's surface carry the vector structure. On flat meshes the curvature is concentrated on the edges, and we average the edge measure obtained from the first variation of area over vertex patches, which gives sparse curvature matrices. At low frequency, four Lagrange-multiplier constraints on the induced charge and on the longitudinal dipole moments of the far field keep the R-ECFOIE accurate in the quasi-static regime, tested on scatterers eight orders of magnitude smaller than the wavelength. Under mesh refinement, the errors decrease at about second order. A Helmholtz resonator, a rounded Fichera cube, and a smoothed Stanford bunny are compared with electric field integral equation (EFIE) solutions, and the first two also with electrostatic references at extreme low frequency.
\end{abstract}


\section{Introduction}
\IEEEPARstart{T}{ime-harmonic} electromagnetic scattering by perfectly electrically conducting (PEC) obstacles is commonly treated with boundary integral equations (BIEs) for the surface current, most notably the electric, magnetic, and combined field integral equations~\cite{volakis2012integral}. Their Galerkin discretization with Rao--Wilton--Glisson (RWG) functions~\cite{rao1982scattering}, commonly referred to as the method of moments, gives rise to dense linear systems whose solution can be accelerated by fast multipole methods or compressed using hierarchical matrices~\cite{bebendorf2000approximation}. The surface-current densities belong to tangential trace spaces with subtle regularity properties, and the corresponding Galerkin theory differs substantially from that of scalar integral equations~\cite{buffa2003galerkin}. Calder\'on preconditioning moreover requires compatible dual basis functions~\cite{buffa2007dual}. The EFIE also suffers from low-frequency breakdown, which can be addressed through loop--star, loop--tree, and quasi-Helmholtz decompositions~\cite{zhao2000integral,merlini2020lowfrequency}. Other low-frequency-stable approaches include current-and-charge~\cite{taskinen2006current} and Debye-source~\cite{epstein2010debye} formulations. High-order Nystr\"om methods have also been developed~\cite{ganesh2006spectrally,bruno2009electromagnetic,hu2021chebyshev,hu2024continuity,hofmann2026continuity}. A different line of work, the so-called field-only surface integral equations~\cite{klaseboer2017nonsingular,sun2020pec,faria2024combined}, avoids surface currents altogether by solving vector Helmholtz problems for the Cartesian components of the field, using a boundary condition on the normal component of the normal derivative of the electric field~\cite{yuffa2018tensorial,faria2024combined}. The Maxwell \`a la Helmholtz framework rigorously establishes that this vector Helmholtz problem is equivalent to the PEC scattering problem and yields indirect~\cite{burbano2026mahelmholtz} and direct~\cite{perez2026direct} second-kind combined-field formulations that are uniquely solvable at all frequencies and have so far been discretized with high-order Nystr\"om methods~\cite{faria2021general}.

The indirect ECFOIE of~\cite{burbano2026mahelmholtz} is developed for $C^{2,\alpha}$-smooth PEC surfaces. It prescribes the tangential electric trace and couples the normal component of its normal derivative to the obstacle's surface mean curvature. The resulting BIE involves only the classical scalar Helmholtz single-layer, double-layer, adjoint double-layer, and hypersingular operators, denoted by $\Sop_k$, $\Kop_k$, $\Kp_k$, and $\Top_k$, respectively, applied componentwise, while tangential and normal projections and the curvature carry the remaining vector structure. We investigate whether it can be accurately discretized with standard continuous P1 densities on flat triangles while retaining the curvature contribution in its variational form. With three components per vertex, a closed mesh with $N$ vertices gives $3N$ unknowns, while an RWG discretization~\cite{rao1982scattering} has one per edge, that is, $3N-6$ on a genus-zero surface by Euler's formula.

The use of standard flat triangular meshes enables practical ECFOIE computations on a broader range of smooth surfaces, including surfaces given only implicitly as level sets, which can be meshed with general-purpose generators such as Gmsh~\cite{geuzaine2009gmsh}. It does not, however, remove the regularity requirement on the underlying geometry, which is a property of the formulation and not of its discretization. Genuinely sharp edges remain outside its scope. Indeed, the componentwise Helmholtz setting requires fields in $H^1$, whereas PEC fields on polyhedral obstacles lie naturally in $H(\mathrm{curl})\cap H(\mathrm{div})$ but generally fail to have $H^1$ regularity because of edge singularities~\cite{costabeldauge2000singularities}. Consequently, nodal $H^1$-conforming discretizations do not in general converge to such fields~\cite{bonnetbendhia1999singular}. Section~\ref{sec:regularity} studies this on surfaces of controlled regularity. The method is thus intended for smooth or rounded scatterers. For these, it is stable at low frequencies without loop--star or quasi-Helmholtz decompositions, uses only scalar Helmholtz operators, so that existing Helmholtz fast multipole and hierarchical-matrix codes can be reused, and gives iteration counts that do not grow under mesh refinement with less memory than an RWG EFIE.

We evaluate the singular integrals with Sauter--Schwab quadrature~\cite{sauter2011bem} and the hypersingular operator in weak form through Maue's identity, discretize the mean curvature variationally following~\cite{meyer2003discrete,sullivan2008curvatures}, impose the low-frequency constraints with Lagrange multipliers, and recover the surface current weakly from the traces. Section~II summarizes the formulation and the constraints, Section~III develops the discretization, Section~IV presents the numerical results, and Section~V discusses scope and limitations.

\section{Electric Field-Only Formulation}
Let $\Omega\subset\R^3$ be a bounded PEC obstacle with a connected boundary $\Gamma$ of class $C^{2,\alpha}$, $0<\alpha<1$, and outward unit normal $\nu$, embedded in a homogeneous medium with permittivity $\epsilon$, permeability $\mu$, and wave impedance $Z=\sqrt{\mu/\epsilon}$. Fields are time-harmonic with time dependence $\e^{-\im\omega t}$, and the wavenumber is $k=\omega\sqrt{\epsilon\mu}$. We use the outgoing Green function
\begin{equation}
 G_k(x,y)=\frac{\e^{\im k|x-y|}}{4\pi|x-y|},\quad x\neq y,\quad k\geq0.
\end{equation}
For a sufficiently smooth scalar density $\phi$ on $\Gamma$, the single- and double-layer potentials are
\begin{equation}\label{eq:layer_potentials}
 \mathcal S_k\phi=\int_\Gamma G_k(\cdot,y)\phi(y)\,\de s_y,\ 
 \mathcal D_k\phi=\int_\Gamma \partial_{\nu_y}G_k(\cdot,y)\phi(y)\,\de s_y
\end{equation}
off $\Gamma$. The boundary operators $\Sop_k$, $\Kop_k$, $\Kp_k$, and $\Top_k$ are the classical boundary integral operators of the Helmholtz equation~\cite{coltonkress1983integral}, with kernels $G_k$, $\partial_{\nu_y}G_k$, $\partial_{\nu_x}G_k$, and $\partial_{\nu_x}\partial_{\nu_y}G_k$, the last in the sense of the Hadamard finite part. We denote by $\gamma^\pm$ and $\partial_\nu^\pm$ the Dirichlet and Neumann traces on~$\Gamma$, that is, the limits of a function and of its normal derivative as~$\Gamma$ is approached from the exterior ($+$) or the interior ($-$), with the outward normal~$\nu$ (understood as a column vector) in both cases. They are applied componentwise to vector fields, and for fields that are smooth across $\Gamma$, such as the incident field, we drop the superscript. They are related to the layer potentials by the jump relations~\cite{coltonkress1983integral}
\begin{equation}\label{eq:jump_relations}
\begin{aligned}
 \gamma^\pm(\mathcal S_k\phi)&=\Sop_k\phi,&
 \gamma^\pm(\mathcal D_k\phi)&=\pm\tfrac12\phi+\Kop_k\phi,\\
 \partial_\nu^\pm(\mathcal S_k\phi)&=\mp\tfrac12\phi+\Kp_k\phi,&
 \partial_\nu^\pm(\mathcal D_k\phi)&=\Top_k\phi.
\end{aligned}
\end{equation}

For the total electric field $E=E^i+E^s$, set $\Pt=\I-\nu\nu^\mathsf T$, $\Pn=\nu\nu^\mathsf T$, and let $\mathscr H=\tfrac12\operatorname{div}_\Gamma\nu$ be the mean curvature (so $\mathscr H=1$ on the unit sphere). The field-only problem in~\cite{burbano2026mahelmholtz} imposes the vector Helmholtz equation $\Delta E^s+k^2E^s=0$ in $\R^3\setminus\overline\Omega$, with the Sommerfeld radiation condition on each component of $E^s$, together with
\begin{equation}\label{eq:pec_traces}
 \Pt\gamma^+E=0,\qquad
 \Pn(\partial_\nu^+E+2\mathscr H\gamma^+E)=0.
\end{equation}
The first boundary condition is the standard PEC condition $\nu\times E=0$. The second replaces the divergence constraint: when $\Pt\gamma^+E=0$, it is equivalent to $\gamma^+(\operatorname{div}E)=0$~\cite[Lemma~3.2]{burbano2026mahelmholtz}. Since $\operatorname{div}E^s$ is then an outgoing Helmholtz solution with zero Dirichlet data, it vanishes, so $E$ is divergence-free and solves Maxwell's equations. For a $C^{2,\alpha}$ boundary, this vector Helmholtz problem is therefore equivalent to the classical PEC scattering problem~\cite[Thm.~3.4]{burbano2026mahelmholtz}.

The scattered field is then represented by the combined-field ansatz
\begin{equation}\label{eq:ansatz}
 E^s=(\mathcal D_k-\im\eta\mathcal S_k)\varphi,
\end{equation}
where $\varphi:\Gamma\to\C^3$ is an auxiliary density and $\eta>0$ is a coupling parameter. (We omit the finite-rank point-source term of~\cite{burbano2026mahelmholtz}, which on flat triangulations only reduces the induced charge, and use the constraints of Section~\ref{sec:constraints} instead.) Imposing~\eqref{eq:pec_traces} on~\eqref{eq:ansatz} gives the ECFOIE system of~\cite{burbano2026mahelmholtz}. To write it in compact form, we introduce the operators
\begin{equation}\label{eq:ops_AB}
 \Aop:=2(\Kop_k-\im\eta\Sop_k),\qquad
 \Bop:=-4(\Top_k-\im\eta\Kp_k).
\end{equation}
By the jump relations~\eqref{eq:jump_relations}, $\gamma^+E^s=\tfrac12(\varphi+\Aop\varphi)$, and we write
\begin{equation}\label{eq:psi_def}
 \psi:=\varphi+\Aop\varphi=2\gamma^+E^s.
\end{equation}
With these definitions, the ECFOIE reads
\begin{equation}\label{eq:ecfoie}
\begin{aligned}
 &\im\eta(\Pt\psi+\Pn\varphi)+2\mathscr H\psi-\tfrac12\Pn\Bop\varphi\\
 &\qquad=-2(\im\eta\Pt+2\mathscr H)\gamma E^i-2\Pn\partial_\nu E^i.
\end{aligned}
\end{equation}
Its normal part contains the hyperersingular operator $\Top_k$ through $\Bop$, so~\eqref{eq:ecfoie} is not of the second kind. The regularized ECFOIE (R-ECFOIE) of~\cite{burbano2026mahelmholtz} applies the static single-layer operator $\Sop_0$, that is, $\Sop_k$ with $k=0$, as a Calder\'on regularizer to the normal part. The resulting second-kind equation is
\begin{equation}\label{eq:recfoie}
\begin{aligned}
 &\Pt\psi+\Sop_0\Pn\bigl[\Bop\varphi-2\im\eta\varphi-4\mathscr H\psi\bigr]\\
 &\qquad=-2\Pt\gamma E^i+4\Sop_0\Pn\bigl[\partial_\nu E^i+2\mathscr H\gamma E^i\bigr].
\end{aligned}
\end{equation}

\subsection{Low Frequencies and the Constrained Formulation}\label{sec:constraints}
Each Cartesian component of~\eqref{eq:ansatz} is a radiating Helmholtz solution, so the ansatz also represents fields that are not divergence-free. The exact density gives a Maxwell field, but a discrete density need not. The far-field pattern $F$ of~\eqref{eq:ansatz}, defined by $E^s(x)\sim\e^{\im k|x|}F(\widehat x)/|x|$ as $|x|\to\infty$, with $\widehat x=x/|x|$ on the unit sphere $\mathbb{S}^2$, is
\begin{equation}\label{eq:farfield}
 F(\widehat x)=\frac{1}{4\pi}\int_\Gamma\e^{-\im k\widehat x\cdot y}\bigl(-\im k\,\widehat x\cdot\nu(y)-\im\eta\bigr)\varphi(y)\,\de s_y.
\end{equation}
A true Maxwell far field is transverse, $\widehat x\cdot F=0$~\cite{coltonkress1983integral}, and $\im k\,\widehat x\cdot F$ is the far field of $\operatorname{div}E^s$.

Three defects of a discrete density dominate as $k\to0$. First, away from the scatterer, the single-layer potential contributes the monopole term $-\im\eta\,\e^{\im k|x|}\int_\Gamma\varphi\,\de s/(4\pi|x|)$. Its coefficient vanishes for the exact density as $k\to0$, but for a discrete density it is of the size of the discretization error and dominates the quasi-static field, which decays like $|x|^{-3}$. Taking $\eta=k$, which mirrors the factor $k$ of the current term of the Stratton--Chu formula $E^s=-\nabla\mathcal S_k[\nu\cdot\gamma^+E]+\im kZ\,\mathcal S_k[J]$~\cite{coltonkress1983integral}, with $J$ the surface current, makes it $O(k)$. Second, with $\eta=k$, constant densities become a near-nullspace: for $k=\eta=0$ they produce no exterior field and are annihilated by the left-hand sides of~\eqref{eq:ecfoie} and~\eqref{eq:recfoie}. Third, the induced charge, proportional to
\begin{equation}\label{eq:charge}
 Q(\varphi)=\int_\Gamma\nu\cdot\gamma^+E^s\,\de s=\tfrac12\int_\Gamma\nu\cdot(\varphi+\Aop\varphi)\,\de s,
\end{equation}
vanishes for the exact solution, since $\operatorname{div}E^i=0$ in $\Omega$ and the conductor is uncharged, but not for a discrete density, and its Coulomb field, proportional to $Q\,x/(4\pi|x|^3)$, decays like $|x|^{-2}$, more slowly than the quasi-static field.

To address these issues, we therefore impose the following constraints:
\begin{equation}\label{eq:constraints}
\begin{aligned}
 b_i(\varphi)&=\int_{\mathbb{S}^2}\widehat x_i\,\widehat x\cdot F(\widehat x)\,\de s_{\widehat x}=0,\quad i\in\{1,2,3\},\\
 b_4(\varphi)&=Q(\varphi)=0.
\end{aligned}
\end{equation}
The first three are the dipole moments of the far field of $\operatorname{div}E^s$. Both conditions hold for the exact solution at every frequency. Since $\int_{\mathbb{S}^2}\widehat x_i\widehat x_j\,\de s=4\pi\delta_{ij}/3$, for $\eta=k$ the first three reduce to $b_i(\varphi)=-\tfrac{\im k}{3}\int_\Gamma\varphi_i\,\de s+O(k^2)$ as $k\to0$. They thus remove both the $|x|^{-1}$ term and the constant near-nullspace, while $b_4$ removes the net charge.

We impose~\eqref{eq:constraints} with Lagrange multipliers, as is standard for constraints in Galerkin boundary element methods, for instance for the total charge in electrostatics~\cite{steinbach2008numerical}. Writing the R-ECFOIE~\eqref{eq:recfoie} as $\mathsf L\varphi=f$, we collect the functionals in $\mathbf b=(b_1,\dots,b_4)$ and seek $\varphi\in H^{1/2}(\Gamma)^3$ and $\boldsymbol\lambda\in\C^4$ such that
\begin{equation}\label{eq:constrained}
\begin{aligned}
 \langle\mathsf L\varphi,v\rangle+\boldsymbol\lambda\cdot\mathbf b(v)&=\langle f,v\rangle\quad\text{for all }v\in H^{1/2}(\Gamma)^3,\\
 \mathbf b(\varphi)&=\mathbf 0,
\end{aligned}
\end{equation}
with the bilinear pairing $\langle\cdot,\cdot\rangle$ of Section~\ref{sec:galerkin}. The exact density solves~\eqref{eq:constrained} with $\boldsymbol\lambda=\mathbf 0$, so the constraints act only on the discrete solution. We take $\eta=k$ at all frequencies. The same construction applies to the ECFOIE~\eqref{eq:ecfoie}, in different function spaces, but with $\eta=k$ it appears to be inaccurate at low frequency (Section~\ref{sec:numerics_low}), so for $kD<1$ we use the R-ECFOIE. The three defects above dominate the error only for $kD\ll1$, with $D$ the diameter of $\Omega$, and for $kD\geq1$ the constraints are in general not needed. They remain consistent there, but they add four unnecessary unknowns to the linear system, so we do not impose them for $kD\geq1$.

\section{Galerkin Discretization}\label{sec:galerkin}
Let $\Gamma_h$ be a closed, consistently oriented, shape-regular polyhedral approximation of $\Gamma$, with flat triangular faces $\mathcal T_h$ and vertices $\{x_a\}_{a=1}^N$. We use $\nu_h$ for its outward unit normal, which is constant on each face. The scalar trial and test space on $\Gamma_h$ is the continuous P1 space $V_h=\operatorname{span}\{\phi_a\}_{a=1}^N\subset C^0(\Gamma_h)$. Each Cartesian component of the vector density is represented in $V_h$, giving $3N$ degrees of freedom. %
Throughout, $\langle u,v\rangle_{\Gamma}=\int_{\Gamma}u\cdot v\,\de s$ denotes the bilinear (unconjugated) pairing on a surface $\Gamma$, with $u\cdot v$ the ordinary product for scalar functions. It is understood as a duality pairing when $u$ is less regular, and as $\langle\mu,v\rangle=\int v\,\de\mu$ when $u=\mu$ is a measure on the surface. The subscript is omitted when the surface is clear. For a scalar boundary operator $\mathsf X\in\{\Aop,\Bop,\Sop_0\}$ constructed on $\Gamma_h$, its corresponding Galerkin matrix $X_h$ has entries
\begin{equation}\label{eq:galerkin}
 (X_h)_{ab}=\langle\mathsf X\phi_b,\phi_a\rangle_{\Gamma_h}.
\end{equation}
This gives the matrices $A_h$, $B_h$, and $S_{0,h}$, the first two each assembled in a single element-pair loop. For $\Bop$, which contains the hypersingular operator $\Top_k$, this is a duality pairing, evaluated through Maue's identity on $\Gamma_h$~\cite{steinbach2008numerical}. The mass matrix, $[M_h]_{ab}=\int_{\Gamma_h}\phi_a\phi_b\,\de s$, is assembled using the exact element formula for affine triangles.

As is standard in modern BEM codes, the singular integrals over touching triangles are evaluated with Sauter--Schwab transformations~\cite{sauter2011bem}, and for the nearly singular integrals over nearby disjoint pairs the part $(1-k^2R^2/2)/(4\pi R)$ of $G_k$, with $R=|x-y|$, and its gradient are integrated analytically over the source triangle~\cite{wilton1984potential}. The near-field criterion and the quadrature orders follow the default BEAST Maxwell assembly~\cite{beastjl}, so that both codes in Section~\ref{sec:numerics} use the same rules. The operator matrices are compressed as hierarchical matrices with HMatrices.jl~\cite{hmatricesjl}, using adaptive cross approximation with tolerance $10^{-6}$ and leaf size 64 unless stated otherwise.

\subsection{Discrete ECFOIE and R-ECFOIE}\label{sec:discrete_equations}
The curvature terms in~\eqref{eq:ecfoie} and~\eqref{eq:recfoie} act on $\psi=\varphi+\Aop\varphi=2\gamma^+E^s$, which we approximate by a P1 vector function $\psi_h\in V_h^3$. Write $\varphi_{i,a}$, $i\in\{1,2,3\}$, $a\in\{1,\dots,N\}$, for the coefficient of the $i$th Cartesian component of $\varphi_h$ at the vertex $x_a$, and let $\boldsymbol\varphi_i=(\varphi_{i,1},\dots,\varphi_{i,N})^\mathsf T\in \mathbb C^{N}$. The matrices $A_h$ and $B_h$ of Section~\ref{sec:galerkin} are mapped to P1 coefficients by mass solves, so that the coefficients of $\psi_h$ are
\begin{equation}\label{eq:discrete_psi}
 \boldsymbol\psi_i:=\boldsymbol\varphi_i+M_h^{-1}A_h\boldsymbol\varphi_i.
\end{equation}

For a scalar $N\times N$ matrix $X_h$, we write $X_h^{(3)}=I_3\otimes X_h$, where $\otimes$ denotes the Kronecker product. Thus, $X_h^{(3)}$ is block diagonal with three copies of $X_h$ and applies $X_h$ to each Cartesian component. Let $P_{\nu,h}$ be the weak facewise discretization of the normal projector $\Pn=\nu\nu^\top$, with entries
\begin{equation}
 [P_{\nu,h}]_{(i,a),(j,b)}=\int_{\Gamma_h}
 \phi_a\phi_b(\nu_h)_i(\nu_h)_j\,\de s.
\end{equation}
The weak tangential projector is $P_{t,h}=M_h^{(3)}-P_{\nu,h}$.
The weak approximation of $\Pn\Bop\varphi$ is
\begin{equation}\label{eq:discrete_normal_B}
 \mathbf r_\nu:=P_{\nu,h}(M_h^{(3)})^{-1}B_h^{(3)}\boldsymbol\varphi.
\end{equation}
Below, the moment vector $\mathbf g\in\C^{3N}$ of a vector field $g$ on $\Gamma_h$ has the entries $[\mathbf g]_{(j,a)}=\int_{\Gamma_h}\phi_a g_j\,\de s$. It contains Galerkin moments rather than nodal values, and $(M_h^{(3)})^{-1}\mathbf g$ gives the coefficients of the $L^2$ projection of $g$ onto $V_h^3$. 

The curvature enters the scheme through two curvature matrices $H_h^{(3)},H_{n,h}\in\C^{3N\times3N}$ and two curvature moment vectors $\mathbf f_{\mathscr H},\mathbf f_{\mathscr H,\nu}\in\C^{3N}$, constructed in Section~\ref{sec:edge_curvature}. The matrices represent in Galerkin form the multiplication operators $u\mapsto\mathscr H u$ and $u\mapsto\Pn\mathscr H u$, respectively, while the vectors contain the Galerkin moments of the corresponding incident-field terms $\mathscr H\gamma E^i$ and $\Pn\mathscr H\gamma E^i$. The discrete ECFOIE, i.e. the Galerkin counterpart of~\eqref{eq:ecfoie}, is then
\begin{equation}\label{eq:discrete_ecfoie}
\begin{aligned}
 \im\eta(P_{t,h}\boldsymbol\psi+P_{\nu,h}\boldsymbol\varphi)
 +2H_h^{(3)}\boldsymbol\psi-\tfrac12\mathbf r_\nu
 ={}&-2\im\eta\,\mathbf f_t-2\mathbf f_\nu\\
 &-4\mathbf f_{\mathscr H},
\end{aligned}
\end{equation}
where $\mathbf f_t$ and $\mathbf f_\nu$ are the moment vectors of $\Pt E^i$ and $\Pn\partial_\nu E^i$, respectively, with the facewise normal $\nu_h$, computed with a high-order triangle rule. For the R-ECFOIE~\eqref{eq:recfoie}, define
\begin{equation}\label{eq:discrete_mu}
 \boldsymbol\mu:=-2\im\eta P_{\nu,h}\boldsymbol\varphi
      +\mathbf r_\nu-4H_{n,h}\boldsymbol\psi.
\end{equation}
We discretize the regularized normal block in the product form obtained directly from~\eqref{eq:recfoie}:
\begin{equation}\label{eq:discrete_product}
\begin{aligned}
 &P_{t,h}\boldsymbol\psi+S_{0,h}^{(3)}(M_h^{(3)})^{-1}\boldsymbol\mu\\
 &\quad=-2\mathbf f_t+4S_{0,h}^{(3)}(M_h^{(3)})^{-1}(\mathbf f_\nu+2\mathbf f_{\mathscr H,\nu}),
\end{aligned}
\end{equation}
with $\mathbf f_{\mathscr H,\nu}$ from~\eqref{eq:curvature_load}. A discretization through the Calder\'on identity $-4\Sop_0\Top_0=\I-4\Kop_0^2$~\cite[Thm.~5.8]{burbano2026mahelmholtz} gives nearly the same results but is up to twice as slow. As in~\cite{bruno2009electromagnetic,boubendir2013wave}, $\Sop_0$ may also be replaced by $\Sop_\kappa$, $\kappa=k+\im k^{1/3}$, which better regularizes $\Top_k$ at high frequencies. Since for $\operatorname{Im}\kappa>0$, $\Sop_\kappa$ is invertible, indeed coercive up to a rotation, and the R-ECFOIE remains of the second kind and uniquely solvable. We use $\Sop_\kappa$ only in Fig.~\ref{fig:sphere_sweep} and Table~\ref{tab:cavity_iters}, since the low-frequency results rely on $\Sop_0$.

\subsection{Discrete Curvature}\label{sec:edge_curvature}
The curvature terms of~\eqref{eq:ecfoie} and~\eqref{eq:recfoie} require $\mathscr H$ on $\Gamma_h$, where it vanishes on each face and is concentrated on the edges. We recover it from the first variation of area~\cite{meyer2003discrete,sullivan2008curvatures}. For a smooth surface and a vector field $V$, with $\mathscr H=\tfrac12\operatorname{div}_\Gamma\nu$, the derivative of the area of $\{x+tV(x)\}$ at $t=0$ is $\int_\Gamma\operatorname{div}_\Gamma V\,\de s=\int_\Gamma 2\mathscr H\nu\cdot V\,\de s$~\cite{coltonkress1983integral}. On $\Gamma_h$, integrating by parts on each face leaves only edge terms, and the discrete counterpart of $2\mathscr H\nu$ is the edge measure~\cite{sullivan2008curvatures}
\begin{figure}[t]
\centering
\includegraphics{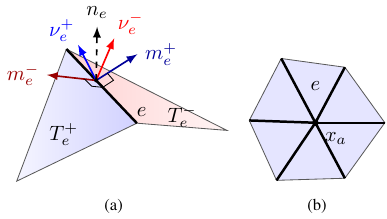}
\caption{(a) Two flat faces $T_e^\pm$ sharing the edge $e$, with outward normals $\nu_e^\pm$, their unit bisector $n_e$, and the outward conormals $m_e^\pm$. (b) Triangles incident to the vertex $x_a$, whose edges $e\ni x_a$ (thick) contribute to the nodal curvature $\mathscr H_a$.}
\label{fig:edge_conormals}
\end{figure}
\begin{equation}\label{eq:curvature_vector}
 \sum_e(m_e^++m_e^-)\,\delta_e=2\sum_e h_en_e\,\delta_e,\quad h_e:=\sin(\theta_e/2),
\end{equation}
where $m_e^\pm$ are the outward unit conormals of the two faces adjacent to the edge $e$ (Fig.~\ref{fig:edge_conormals}(a)), $n_e$ is the unit bisector of their normals, $\theta_e$ is the fold angle, positive at convex edges, and $\delta_e$ is the arclength measure on $e$, so that $\langle\delta_e,v\rangle=\int_e v\,\de\ell$. The scalar curvature measure is then given by $\mathscr H_h=\sum_e h_e\delta_e$.

Concentrated on the edges, $\mathscr H_h$ would make the curvature terms act only on the edge traces of $\psi_h\in V_h^3$. We instead average it over the vertex patches, which turns the singular edge measure into a P1 field compatible with the trial space and intended to recover the smooth curvature under refinement. Let $A_a=\int_{\Gamma_h}\phi_a\,\de s$ be the area assigned to $x_a$, a third of the area of the triangles incident to $x_a$ (Fig.~\ref{fig:edge_conormals}(b)). The nodal curvature and the P1 curvature field are then defined by
\begin{equation}\label{eq:nodal_curvature}
 \begin{aligned}\mathscr H_a:=&~\frac{\langle\mathscr H_h,\phi_a\rangle}{A_a}
 =\frac{1}{2A_a}\sum_{e\ni x_a}h_eL_e\quad\text{and}\\ 
 \mathscr H_h^{(1)}:=&~\sum_{a=1}^N\mathscr H_a\phi_a\in V_h,
 \end{aligned}
\end{equation}
respectively, where the first sum runs over the edges $e$ incident to $x_a$, $L_e$ is the length of $e$, and the factor $\tfrac12$ comes from $\langle\delta_e,\phi_a\rangle=L_e/2$. Since the P1 basis functions satisfy $\sum_a\phi_a=1$, the total curvature is unchanged, $\int_{\Gamma_h}\mathscr H_h^{(1)}\de s=\langle\mathscr H_h,1\rangle$. We then get that the curvature matrices are given by
\begin{equation}\label{eq:curvature_mass}
\begin{aligned}
 [H_h]_{ab}&=\int_{\Gamma_h}\mathscr H_h^{(1)}\phi_a\phi_b\,\de s,\\
 [H_{n,h}]_{(i,a),(j,b)}&=\int_{\Gamma_h}\mathscr H_h^{(1)}\phi_a\phi_b(\nu_h)_i(\nu_h)_j\,\de s,
\end{aligned}
\end{equation}
and $H_h^{(3)}$ acts on all three components. The integrands above are cubic polynomials on each face and are integrated exactly. The matrices $H_h^{(3)}$ and $H_{n,h}$ discretize the multiplication operators $\psi\mapsto\mathscr H\psi$ and $\psi\mapsto\Pn\mathscr H\psi$, the latter with the facewise normal of the projector $P_{\nu,h}$. The curvature term $2\mathscr H\psi$ of the ECFOIE~\eqref{eq:ecfoie} thus becomes $2H_h^{(3)}\boldsymbol\psi$, and the normal-block term $-4\Pn\mathscr H\psi$ of the R-ECFOIE~\eqref{eq:recfoie} becomes $-4H_{n,h}\boldsymbol\psi$. The incident field supplies the corresponding moment vectors
\begin{equation}\label{eq:curvature_load}
\begin{aligned}
 [\mathbf f_{\mathscr H}]_{(j,a)}&=\int_{\Gamma_h}\mathscr H_h^{(1)}\phi_a E^i_j\,\de s,\\
 [\mathbf f_{\mathscr H,\nu}]_{(j,a)}&=\int_{\Gamma_h}\mathscr H_h^{(1)}\phi_a(\nu_h)_j\bigl(\nu_h\cdot E^i\bigr)\,\de s,
\end{aligned}
\end{equation}
evaluated with the same face rule as the other incident-field moment vectors. The first enters the ECFOIE right-hand side as $-4\mathbf f_{\mathscr H}$, and the second contributes $8S_{0,h}^{(3)}(M_h^{(3)})^{-1}\mathbf f_{\mathscr H,\nu}$ to the regularized right-hand side. %

When the discrete normals converge, this edge measure converges in a weak sense to $2\mathscr H\nu$~\cite{hildebrandt2006convergence}. %

\subsection{Implementation of the Low-Frequency Constraints}\label{sec:constraint_impl}
The functionals~\eqref{eq:constraints} are evaluated on the basis of $V_h^3$, which gives a matrix $\Lambda_h\in\C^{4\times3N}$ with entries $[\Lambda_h]_{j,(i,a)}=b_j(\phi_a\mathbf e_i)$. Each entry $b_j(\phi_a\mathbf e_i)$ of a far-field row is a double integral, over $\Gamma_h$ in~\eqref{eq:farfield} and over $\mathbb{S}^2$ in~\eqref{eq:constraints}. Over $\Gamma_h$ we use the facewise normal $\nu_h$ and the same triangle rule as for the incident-field moment vectors in~\eqref{eq:discrete_ecfoie}. Over $\mathbb{S}^2$ we use a standard product rule of sufficiently high degree. For the charge row, the discrete counterpart of~\eqref{eq:charge} is $Q_h=\tfrac12\mathbf n_h^\mathsf T\boldsymbol\psi=\mathbf q_h^\mathsf T\boldsymbol\varphi$, with $\mathbf n_h$ the moment vector of $\nu_h$. Inserting~\eqref{eq:discrete_psi} and using the symmetry of $M_h$ gives the row vector, computed once,
\begin{equation}\label{eq:charge_row}
 \mathbf q_h=\tfrac12\bigl(\mathbf n_h+(A_h^\mathsf T)^{(3)}(M_h^{(3)})^{-1}\mathbf n_h\bigr),
\end{equation}
which is the fourth row of $\Lambda_h$, that is, $[\Lambda_h]_{4,(i,a)}=[\mathbf q_h]_{(i,a)}$.

Let $L_h\boldsymbol\varphi=\mathbf f_h$ denote the discrete R-ECFOIE~\eqref{eq:discrete_product}, with $L_h$ the matrix of its left-hand side, in which $\boldsymbol\psi$ and $\boldsymbol\mu$ depend linearly on $\boldsymbol\varphi$, and $\mathbf f_h=-2\mathbf f_t+4S_{0,h}^{(3)}(M_h^{(3)})^{-1}(\mathbf f_\nu+2\mathbf f_{\mathscr H,\nu})$. The Galerkin counterpart of~\eqref{eq:constrained} is the saddle-point system~\cite{steinbach2008numerical}
\begin{equation}\label{eq:saddle}
 \begin{bmatrix}L_h&\Lambda_h^\mathsf T\\ \Lambda_h&0\end{bmatrix}
 \begin{bmatrix}\boldsymbol\varphi\\ \boldsymbol\lambda\end{bmatrix}
 =\begin{bmatrix}\mathbf f_h\\ \mathbf 0\end{bmatrix},
\end{equation}
in which $\Lambda_h^\mathsf T$ is the transpose, not the conjugate transpose, because the pairing is bilinear. As for the unconstrained equation, we precondition from the left with $\operatorname{blockdiag}((M_h^{(3)})^{-1},I_4)$. We also scale the constraints symmetrically, replacing $\Lambda_h$ by $\Sigma\Lambda_h$ and $\Lambda_h^\mathsf T$ by $\Lambda_h^\mathsf T\Sigma$, with $\Sigma=\operatorname{diag}(\sigma_1,\dots,\sigma_4)$ and $\sigma_j=\sqrt N/\|(M_h^{(3)})^{-1}\Lambda_h^\mathsf T\mathbf e_j\|_2$, so that the columns of $(M_h^{(3)})^{-1}\Lambda_h^\mathsf T\Sigma$ have norm $\sqrt N$. As noted after~\eqref{eq:constraints}, the far-field rows are $O(k)$ as $k\to0$, and $\Sigma$ removes this factor. The preconditioned matrix is thus the second-kind operator $(M_h^{(3)})^{-1}L_h$ augmented by four rows and columns that remain of comparable size as $k\to0$, so its conditioning does not deteriorate. The system is solved by GMRES. Each iteration costs one application of $L_h$ and $O(N)$ operations for the four constraints.

\subsection{Surface Current}\label{sec:current_recovery}
In many applications, the surface current is itself a quantity of interest. The density $\varphi$ is auxiliary, but the physical surface current $J=\nu\times H^+=(\im kZ)^{-1}\nu\times\operatorname{curl}E^+$ can be recovered weakly from the traces of the solution. With $E_\nu=\nu\cdot\gamma^+E=\nu\cdot(\gamma E^i+\tfrac12\psi)$ and $g=\partial_\nu^+E=\partial_\nu E^i-\tfrac14\Bop\varphi+\tfrac{\im\eta}{2}\varphi$, the identity $\nu\times\operatorname{curl}E^+=\nabla_\Gamma E_\nu-\Pt g$~\cite[Lemma~3.2]{burbano2026mahelmholtz} and integration by parts give, for tangential $w$,
\begin{equation}\label{eq:current_weak_identity}
 \langle J,w\rangle=\frac{1}{\im kZ}\bigl[-\langle E_\nu,\operatorname{div}_\Gamma w\rangle-\langle g,w\rangle\bigr].
\end{equation}
Testing with lowest-order RWG functions gives $J_h$ from a sparse mass-matrix system, in which $g$ is replaced by its P1 approximation $g_h\in V_h^3$, obtained by one mass solve from $B_h^{(3)}\boldsymbol\varphi$ and the moment vector of $\partial_\nu E^i$. The recovered current also permits a field evaluation through the Stratton--Chu formula~\cite{coltonkress1983integral}, but with the constraints of Section~\ref{sec:constraints} this is not needed at low frequency, and we evaluate all fields with~\eqref{eq:ansatz}.

\section{Numerical Results}\label{sec:numerics}
We consider the unit sphere, a Helmholtz resonator, surfaces of controlled regularity, a rounded Fichera cube, and a smoothed Stanford bunny, in order of increasing geometric difficulty. Unless stated otherwise, we solve the R-ECFOIE in the product form~\eqref{eq:discrete_product} with $k=\eta=\pi$, for the unit plane wave $E^i(x)=\mathbf e_1\,\e^{\im kx_3}$, polarized along $x$ and propagating along $+z$. For $kD<1$ we add the constraints~\eqref{eq:constraints} and solve the saddle-point system~\eqref{eq:saddle}. In all low-frequency tests we take $D=10^{-8}\lambda$, with $\lambda=2\pi/k$, where the EFIE suffers from low-frequency breakdown~\cite{zhao2000integral}. As $kD\to0$, the scattered field converges, with relative error of order $kD$, to the electrostatic field of the uncharged conductor in the uniform field $E^i(0)$~\cite{dassios2000low}, which we use as the reference. For the unit sphere, this is the dipole field $[3\widehat x(\widehat x\cdot E^i(0))-E^i(0)]/|x|^3$. For other geometries, we compute it from the Galerkin Laplace single-layer equation with zero total charge~\cite{steinbach2008numerical}. Errors are evaluated at 100 nearly uniformly distributed directions $\widehat x_j\in\mathbb S^2$, from a Fibonacci lattice: the scattered field $E^s$ at the points $d\,\widehat x_j$, with $d=D$, and the far-field pattern $F$ in~\eqref{eq:farfield} in the directions $\widehat x_j$. With $\|u\|=(\sum_j|u_j|^2)^{1/2}$ and $F_\perp=(I-\widehat x\widehat x^\mathsf T)F$, we report
\begin{equation}\label{eq:farfield_metrics}
 \varepsilon_E=\frac{\|E^s-E^{s,\mathrm{ref}}\|}{\|E^{s,\mathrm{ref}}\|},\ \ 
 \varepsilon_F=\frac{\|F_\perp-F_\perp^{\mathrm{ref}}\|}{\|F_\perp^{\mathrm{ref}}\|},\ \ 
 \varepsilon_L=\frac{\|\widehat x\cdot F\|}{\|F\|},
\end{equation}
where the superscript ``ref'' denotes a reference solution and~$\varepsilon_L$, which needs none, measures the violation of far-field transversality. The references are the Mie series for the sphere at $D=\lambda$, the electrostatic field described above at low frequency, and, for the other geometries at $D=\lambda$ and above, BEAST.jl EFIE solutions with RWG functions~\cite{beastjl}.

\subsection{PEC Sphere}\label{sec:numerics_sphere}
We first consider scattering of the default plane wave by the unit sphere with $D=\lambda$, that is, $k=\eta=\pi$. The reference is the Mie-series solution, with $d=2$, and the meshes are listed in Table~\ref{tab:sphere}. We solve the ECFOIE and the R-ECFOIE, each without and with the constraints~\eqref{eq:constraints}, using the same assembled operators on each mesh, so that the differences are due to the formulations alone. Both are preconditioned by the inverse mass matrix and solved with unrestarted GMRES, with relative tolerance $10^{-6}$ and at most 400 iterations.

\begin{table}[h!]
\caption{Sphere with $D=\lambda$: Near-Field Error $\varepsilon_E$ (\%) and GMRES Iterations, Without and With the Constraints}
\label{tab:sphere}
\centering
\small
\setlength{\tabcolsep}{3pt}
\begin{tabular}{crrrrr}
\toprule
 & & \multicolumn{2}{c}{R-ECFOIE} & \multicolumn{2}{c}{ECFOIE} \\
$h$ & $N$ & without & with & without & with \\
\midrule
0.250 & 271  & 3.47 (24)  & 3.43 (30)  & 2.78 (61)   & 2.76 (74)   \\
0.175 & 508  & 1.73 (24)  & 1.68 (29)  & 1.44 (69)   & 1.44 (85)   \\
0.125 & 1050 & 0.77 (24)  & 0.75 (28)  & 0.67 (75)   & 0.66 (92)   \\
0.100 & 1578 & 0.50 (24)  & 0.49 (28)  & 0.44 (80)   & 0.44 (99)   \\
0.070 & 3115 & 0.25 (23)  & 0.24 (27)  & 0.22 (93)   & 0.22 (106)  \\
0.050 & 6064 & 0.12 (22)  & 0.12 (26)  & 0.11 (107)  & 0.11 (117)  \\
0.040 & 9514 & 0.077 (22) & 0.074 (26) & 0.072 (119) & 0.071 (130) \\
\bottomrule
\end{tabular}
\end{table}

Table~\ref{tab:sphere} shows that the field errors decrease at second order in $h$, that is, like $N^{-1}$ (fitted slopes $2.00$--$2.09$ for $h\le0.175$). This rate is consistent with the $O(h^2)$ errors of both P1 interpolation and the flat-triangle approximation of~$\Gamma$. The two equations reach nearly the same accuracy, but the R-ECFOIE iteration count stays bounded under refinement, whereas that of the ECFOIE grows steadily. The constraints change the errors only slightly and add a few iterations, so at this frequency they are not needed, as noted in Section~\ref{sec:constraints}. Without the curvature term (not shown), the errors of both equations stagnate at about $70\%$ on all meshes.

To give an idea of the overall cost of the proposed formulations, we compare with the EFIE of BEAST.jl~\cite{beastjl} on the same sphere, without the constraints, on meshes with $N=3115$ to $30966$ vertices\label{sec:beast}. Our formulations then have $3N$ unknowns and the EFIE has $3N-6$, one per edge. Both codes use the quadrature of Section~\ref{sec:galerkin} and unrestarted GMRES with tolerance $10^{-6}$, Jacobi-preconditioned for the EFIE. All matrices are compressed, with HMatrices.jl for our formulations and FastBEAST.jl for the EFIE, and are also assembled densely up to $N=9514$. Each code uses the loosest compression that changes $\varepsilon_E$ by less than $1\%$: ACA tolerance $10^{-3}$ and leaf size $128$ for HMatrices.jl, and tolerance $10^{-4}$, minimum box size $200$, and maximum rank $50$ for FastBEAST.jl. Times are medians of three runs on an Apple M3 Max with $10$ threads.

\begin{table}[h!]
\caption{Sphere with $D=\lambda$: Comparison with the BEAST EFIE; Times in Seconds}
\label{tab:beast}
\centering
\footnotesize
\setlength{\tabcolsep}{3.5pt}
\begin{tabular}{rlrrrrrr}
\toprule
 & & & \multicolumn{2}{c}{H-matrix} & \multicolumn{2}{c}{Dense} & \\
\cmidrule(lr){4-5}\cmidrule(lr){6-7}
$N$ & Method & It. & Assem. & Solve & Assem. & Solve & $\varepsilon_E$ (\%) \\
\midrule
3115 & ECFOIE   & 93  & 5.1  & 0.8  & 3.6  & 0.9  & 0.222 \\
     & EFIE     & 240 & 5.0  & 2.7  & 4.4  & 3.5  & 0.199 \\
     & R-ECFOIE & 23  & 8.7  & 0.4  & 7.4  & 0.3  & 0.247 \\
\midrule
6064 & ECFOIE   & 107 & 13.6 & 2.6  & 13.2 & 3.0  & 0.113 \\
     & EFIE     & 300 & 12.5 & 8.6  & 17.6 & 15.8 & 0.100 \\
     & R-ECFOIE & 23  & 21.1 & 0.9  & 28.0 & 1.0  & 0.121 \\
\midrule
9514 & ECFOIE   & 119 & 26.8 & 4.8  & 32.7 & 7.2  & 0.0721 \\
     & EFIE     & 353 & 20.4 & 19.9 & 45.6 & 43.1 & 0.0634 \\
     & R-ECFOIE & 22  & 39.5 & 1.4  & 71.0 & 2.4  & 0.0763 \\
\midrule
30966 & ECFOIE   & 162 & 123.1 & 19.5  & --- & --- & 0.0228 \\
      & EFIE     & 554 & 108.4 & 142.9 & --- & --- & 0.0192 \\
      & R-ECFOIE & 22  & 182.7 & 4.9   & --- & --- & 0.0224 \\
\bottomrule
\end{tabular}
\end{table}

Table~\ref{tab:beast} shows that the three discretizations reach nearly the same accuracy on each mesh, with the EFIE slightly ahead, and that compression changes the errors only marginally. The iteration counts differ markedly. As expected, that of the EFIE grows steadily under refinement, that of the ECFOIE grows more slowly, and that of the R-ECFOIE stays constant, so the R-ECFOIE solve is the fastest on every mesh and its advantage grows with $N$. Our assembly is more expensive, because it involves several scalar operators and, for the R-ECFOIE, the additional static single-layer operator of~\eqref{eq:discrete_product}. Even so, with compression the ECFOIE has the lowest total time on every mesh, and the R-ECFOIE overtakes the EFIE as $N$ grows. Our formulations also need less memory: in dense form, the two scalar matrices $A_h$ and $B_h$ occupy less than a quarter of the memory of the EFIE matrix, and with compression the peak memory on the finest mesh is about $30\%$ lower. These timings depend on the implementations and the tuning, but the iteration counts do not. The Jacobi-preconditioned EFIE is the most common RWG formulation. Since the EFIE is of the first kind, its iteration count grows under mesh refinement, so it is the natural reference for the unregularized ECFOIE. The natural reference for the R-ECFOIE is instead the Calder\'on-preconditioned EFIE~\cite{andriulli2008multiplicative}, whose iteration count is likewise bounded under refinement. We compare the two over frequency below.

\begin{figure}[t]
\centering
\includegraphics[width=\columnwidth]{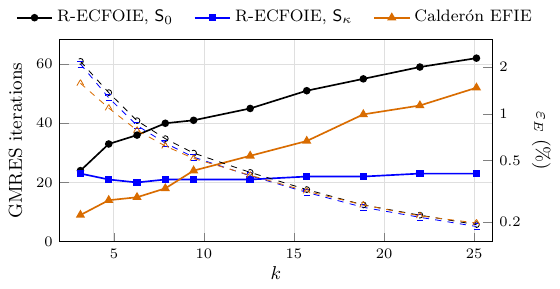}
\caption{Unit sphere, ten points per wavelength: GMRES iterations (solid, left axis) and $\varepsilon_E$ (dashed, right axis) versus $k$ of the R-ECFOIE with the regularizers $\Sop_0$ and $\Sop_\kappa$, $\kappa=k+\im k^{1/3}$, and of the Calder\'on-preconditioned RWG EFIE.}
\label{fig:sphere_sweep}
\end{figure}

To study the dependence on the frequency, we solve the R-ECFOIE with $\Sop_0$ and with $\Sop_\kappa$ for $k=\eta$ from $\pi$ to $8\pi$, that is, $D$ from $\lambda$ to $8\lambda$, on meshes with ten points per wavelength, $h=\lambda/10$, and up to $24082$ vertices. Fig.~\ref{fig:sphere_sweep} shows that with $\Sop_0$ the iteration count grows from $24$ to $62$, whereas with $\Sop_\kappa$ it stays between $20$ and $23$. The errors are nearly the same for both regularizers, and they decrease with $k$ because the meshes also resolve the sphere better. For comparison, we solve the Calder\'on-preconditioned RWG EFIE~\cite{andriulli2008multiplicative} with Buffa--Christiansen dual functions~\cite{buffa2007dual} in BEAST.jl on the same meshes. At $k\le2.5\pi$, it needs fewer iterations and is slightly more accurate. Its iteration count, however, grows from $9$ to $52$ and exceeds that of the R-ECFOIE with $\Sop_\kappa$ from $k=3\pi$ on, where the accuracy of the two is nearly the same.

Next, to study the low-frequency behavior of our formulations, we take $D=10^{-8}\lambda$ for the same sphere\label{sec:numerics_low}, that is, $k=\eta=\pi\cdot10^{-8}$, and solve the R-ECFOIE without and with the constraints~\eqref{eq:constraints}, with the rest of the setup unchanged. The reference is the static dipole field described at the beginning of this section. Table~\ref{tab:low_frequency} reports $\varepsilon_E$ at the radii $d$ and $4d$, the GMRES iterations, and the induced charge $q=|Q_h|/|\Gamma_h|$, with $Q_h$ the discrete counterpart of~\eqref{eq:charge}.
\begin{table}[h!]
\caption{Sphere with $D=10^{-8}\lambda$: R-ECFOIE Near-Field Error $\varepsilon_E$ (\%) at Radii $d$ and $4d$, GMRES Iterations, and Induced Charge $q=|Q_h|/|\Gamma_h|$, Without and With the Constraints}
\label{tab:low_frequency}
\centering
\footnotesize
\setlength{\tabcolsep}{2.5pt}
\begin{tabular}{rrrrcrrrc}
\toprule
 & \multicolumn{4}{c}{Without constraints} & \multicolumn{4}{c}{With constraints} \\
\cmidrule(lr){2-5}\cmidrule(lr){6-9}
$N$ & $d$ & $4d$ & It. & $q$ & $d$ & $4d$ & It. & $q$ \\
\midrule
83   & 343.3  & 1369 & 41 & $2{\cdot}10^{0}$  & 10.00 & 10.34 & 17 & $3{\cdot}10^{-6}$ \\
271  & 123.6  & 453.6 & 40 & $8{\cdot}10^{-1}$ & 2.889 & 2.926 & 14 & $2{\cdot}10^{-6}$ \\
1050 & 219.8  & 869.8 & 42 & $1{\cdot}10^{0}$  & 0.706 & 0.708 & 12 & $2{\cdot}10^{-7}$ \\
3115 & 1465   & 5862  & 44 & $1{\cdot}10^{1}$  & 0.230 & 0.231 & 12 & $6{\cdot}10^{-8}$ \\
9514 & 70.79  & 280.9 & 45 & $5{\cdot}10^{-1}$ & 0.073 & 0.073 & 11 & $2{\cdot}10^{-8}$ \\
\bottomrule
\end{tabular}
\end{table}

Without the constraints, the R-ECFOIE fails at this frequency. The errors are comparatively large, do not decrease under refinement, and grow with the distance to the scatterer, and the induced charge is of the order of the incident field. This is the behavior anticipated in Section~\ref{sec:constraints}. In fact, with $\eta=k$, the constant densities form a near-nullspace, the charge is not controlled, and its Coulomb field, proportional to $Q_h\,x/(4\pi|x|^3)$, decays more slowly than the exact dipole field. With the constraints, the charge vanishes to the solver tolerance, the errors decrease at second order and no longer depend on the distance, and GMRES needs fewer iterations, with a count that does not grow under refinement. With $\eta=k$, the errors of the unregularized ECFOIE remain large with or without the constraints, likely because its first-kind normal part loses the term $\im\eta\Pn\varphi$ as $\eta\to0$, so we use the R-ECFOIE at low frequency.

\subsection{Helmholtz Resonator}\label{sec:cavity}
\begin{figure}[h]
\centering
\includegraphics[width=\columnwidth]{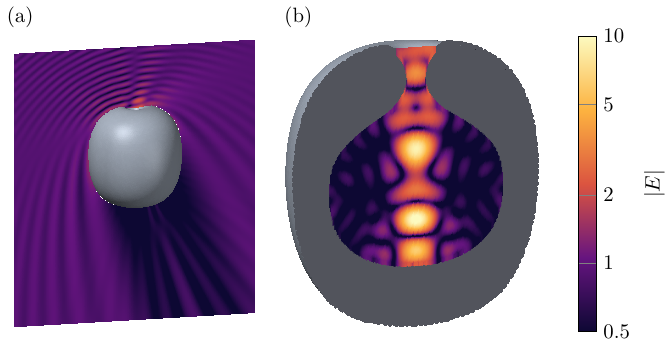}
\caption{Total field $|E|$ on the plane $y=0$ for the Helmholtz resonator at $k=19.5$, where $D\approx8\lambda$, for a plane wave incident at $20^\circ$ from the $-z$ axis in the plane $y=0$, with $E^i$ in that plane. (a) Around the closed body, which produces a shadow below it. (b) Inside the cavity, where the gray band is the cut through the PEC wall.}%
\label{fig:cavity}
\end{figure}

We next consider a smooth axisymmetric body with an open internal chamber, seen in cross-section in Fig.~\ref{fig:cavity}, obtained by revolving a smooth profile curve, a cubic B-spline, about the $z$ axis. Its outer radius is $1.14$, its height $2.45$, and $D\approx2.56$, and a neck of radius $0.121$ connects the mouth, at the top, to a chamber of radius $0.86$, so that the boundary is a single smooth closed surface. Unless stated otherwise, the incident plane wave in this section is $E^i(x)=\mathbf e_1\,\e^{-\im kx_3}$, which enters the mouth along $-z$. The five meshes, with $4660$ to $61860$ vertices, are generated by Gmsh with a curvature-based size field.

We first take $D=\lambda$ and solve both equations without the constraints. The reference is a BEAST EFIE solution, compressed by FastBEAST.jl and solved with tolerance $10^{-8}$, on the 31773-vertex mesh, and $\varepsilon_E$ is evaluated with $d=D$.

\begin{table}[h!]
\caption{Helmholtz Resonator with $D=\lambda$: Near- and Far-Field Errors and Longitudinal Fraction (\%) and GMRES Iterations Against the BEAST EFIE}
\label{tab:cavity}
\centering
\footnotesize
\setlength{\tabcolsep}{2.5pt}
\begin{tabular}{rrrrrrrrr}
\toprule
 & \multicolumn{4}{c}{R-ECFOIE} & \multicolumn{4}{c}{ECFOIE} \\
\cmidrule(lr){2-5}\cmidrule(lr){6-9}
$N$ & $\varepsilon_E$ & $\varepsilon_F$ & $\varepsilon_L$ & It. & $\varepsilon_E$ & $\varepsilon_F$ & $\varepsilon_L$ & It. \\
\midrule
4660            & 0.430 & 0.406 & 0.12   & 37 & 0.310 & 0.285 & 0.11   & 305 \\
6992            & 0.229 & 0.215 & 0.068  & 37 & 0.173 & 0.157 & 0.069  & 302 \\
15642           & 0.084 & 0.077 & 0.029  & 37 & 0.066 & 0.057 & 0.032  & 369 \\
31773$^\ast$    & 0.029 & 0.025 & 0.014  & 36 & 0.030 & 0.025 & 0.016  & 442 \\
61860           & 0.018 & 0.016 & 0.0069 & 35 & 0.026 & 0.024 & 0.0088 & 527 \\
\bottomrule
\multicolumn{9}{l}{\footnotesize $^\ast$Mesh of the reference.}
\end{tabular}
\end{table}

Table~\ref{tab:cavity} shows the same behavior as for the sphere. Both equations converge at about second order until the errors reach the level of the reference, on the finest meshes, and the longitudinal fraction $\varepsilon_L$ decreases in the same way. The R-ECFOIE iteration count is again independent of the mesh, whereas that of the ECFOIE grows under refinement and is already an order of magnitude larger, as is that of the EFIE.

\begin{table}[h]
\caption{Helmholtz Resonator with $D=10^{-8}\lambda$: R-ECFOIE Near-Field Error $\varepsilon_E$ (\%) at Radii $D$ and $4D$, GMRES Iterations, and Induced Charge $q$, Without and With the Constraints}
\label{tab:cavity_low}
\centering
\footnotesize
\setlength{\tabcolsep}{2.5pt}
\begin{tabular}{rrrrcrrrc}
\toprule
 & \multicolumn{4}{c}{Without constraints} & \multicolumn{4}{c}{With constraints} \\
\cmidrule(lr){2-5}\cmidrule(lr){6-9}
$N$ & $D$ & $4D$ & It. & $q$ & $D$ & $4D$ & It. & $q$ \\
\midrule
4660  & 15.92 & 63.08 & 75 & $7{\cdot}10^{-2}$ & 0.377 & 0.495 & 28 & $1{\cdot}10^{-7}$ \\
6992  & 1.679 & 3.366 & 77 & $3{\cdot}10^{-3}$ & 0.206 & 0.327 & 28 & $5{\cdot}10^{-8}$ \\
15642 & 22.49 & 90.09 & 64 & $1{\cdot}10^{-1}$ & 0.082 & 0.136 & 28 & $7{\cdot}10^{-8}$ \\
31773 & 123.7 & 496.1 & 82 & $5{\cdot}10^{-1}$ & 0.033 & 0.066 & 28 & $6{\cdot}10^{-8}$ \\
\bottomrule
\end{tabular}
\end{table}

At $D=10^{-8}\lambda$, the reference is the electrostatic field computed on the 61860-vertex mesh. Table~\ref{tab:cavity_low} confirms the observations made for the sphere, i.e., without the constraints the induced charge is not controlled and the errors are large and erratic, whereas with them the charge vanishes, the errors decrease at second order, and the iteration count is the same on all meshes.

At higher frequencies, the neck acts as a circular waveguide and above $k\approx 15.2$ the planewave enters the chamber. We locate resonances by scanning $k$ and monitoring the mean of $|E|^2$ at points inside the chamber, using the 15642-vertex mesh. Fig.~\ref{fig:cavity} shows the total field at the resonance $k=19.5$ on the plane $y=0$, for incidence at $20^\circ$ from the $-z$ axis in the plane $y=0$, with $E^i$ in that plane. Inside, Fig.~\ref{fig:cavity}(b), a standing wave fills the neck and the chamber, the oblique incidence breaks its axial symmetry, and the field reaches about $10$ times the incident amplitude. Outside, Fig.~\ref{fig:cavity}(a), the incident and reflected waves interfere at the cavity opening and the body produces a shadow below it. At these frequencies we do not attempt to produce an EFIE reference field, so we assess the solution with two quantities that vanish for an exact Maxwell field: the far-field longitudinal fraction $\varepsilon_L$ and the divergence ratio $\|\nabla\cdot E^s_h\|/(k\|E^s_h\|)$ at the points inside the chamber. Both stay below $3\times10^{-3}$ over the frequency range $15.5\le k\le20$.

\begin{table}[h]
\caption{Helmholtz Resonator, 15642 Vertices: GMRES Iterations of the R-ECFOIE and of the Scalar Combined-Field Equation~\eqref{eq:scalar_cfie}; Resonances Marked $^\ast$; Incidence Along $-z$, Except at $20^\circ$ in the Last Column}
\label{tab:cavity_iters}
\centering
\footnotesize
\setlength{\tabcolsep}{3pt}
\begin{tabular}{lrrrrrrr}
\toprule
$k$ & 15.25 & 15.5$^\ast$ & 18.75$^\ast$ & 19.25 & 19.5$^\ast$ & 20 & 19.5$^\ast$($20^\circ$) \\
\midrule
R-ECFOIE, $\Sop_0$ & 323 & 542 & 773 & 647 & 750 & 943 & 807 \\
R-ECFOIE, $\Sop_\kappa$ & 168 & 402 & 610 & 402 & 500 & 617 & 631 \\
Scalar CFIE  & 54  & 57  & 143 & 138 & 128 & 143 & 149 \\
\bottomrule
\end{tabular}
\end{table}

To put the iteration counts in perspective, we also solve the scalar sound-soft problem for the plane wave $u^i=\e^{\im k\,\widehat q\cdot x}$ with $\widehat q=-\mathbf e_3$, as for $E^i$, or tilted by $20^\circ$ in the plane $y=0$ for the last column of Table~\ref{tab:cavity_iters}, using the combined-field equation
\begin{equation}\label{eq:scalar_cfie}
\big(\tfrac12\I+\Kop_k-\im\eta\Sop_k\big)\sigma=-u^i,\qquad \eta=k,
\end{equation}
for the density of $u^s=(\mathcal D_k-\im\eta\mathcal S_k)\sigma$, discretized with the same P1 elements, mass-matrix preconditioner, and GMRES tolerance. Table~\ref{tab:cavity_iters} shows that the resonator is hard already for this scalar problem, whose iteration count increases markedly once $k$ is well above the cutoff, as expected for boundary integral equations on geometries that trap waves~\cite{marchand2022applying}. The R-ECFOIE needs several times more iterations, and its count is higher at some of the electromagnetic resonances, which need not coincide with those of the scalar problem. With $\Sop_\kappa$, the R-ECFOIE needs $21$--$48\%$ fewer iterations, and its far field differs from that for $\Sop_0$ by $\varepsilon_F\le0.4\%$ away from the resonances and by $0.8$--$5.7\%$ at them. Resonances are sensitive to discretization errors, and the mesh has only about seven points per wavelength there.

\subsection{Surfaces of Controlled Regularity}\label{sec:regularity}
\begin{figure}[t]
\centering
\includegraphics[width=\columnwidth]{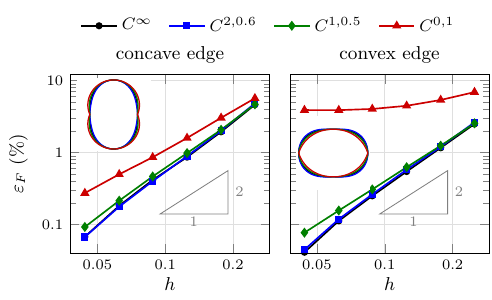}
\caption{Surfaces $r(\omega)\,\omega$ with $r=1+0.45|\omega_3|^p$ (left, concave edge for $p=1$) and $r=1.45-0.45|\omega_3|^p$ (right, convex edge for $p=1$) at $k=\eta=\pi$: far-field error $\varepsilon_F$ of the R-ECFOIE versus the mesh size $h$, for $p=2$, $2.6$, $1.5$, and $1$, which give the regularity in the legend. The insets show the cross-sections of the surfaces. The reference error is about $0.04\%$.}
\label{fig:regularity}
\end{figure}

To see how the regularity of $\Gamma$ affects the convergence, we consider the surfaces $\{r(\omega)\,\omega:\omega\in\mathbb S^2\}$ with $r=1+0.45|\omega_3|^p$ or $r=1.45-0.45|\omega_3|^p$, both with diameter $D=2.9$. Their regularity at the equator $\omega_3=0$ is controlled by $p$: $C^\infty$ for $p=2$, $C^{2,0.6}$ for $p=2.6$, $C^{1,0.5}$ for $p=1.5$, and Lipschitz for $p=1$, with a mild edge along the equator. The edge is concave for the first family and convex for the second (insets of Fig.~\ref{fig:regularity}). Near an edge with exterior opening angle $\Omega$, the field behaves like $d^{\pi/\Omega-1}$, with $d$ the distance to the edge~\cite{costabeldauge2000singularities}. It is therefore in $H^1$ at the concave edge, with $\Omega\approx132^\circ$, but not at the convex one, with $\Omega\approx215^\circ$. We solve for the default plane wave with $k=\eta=\pi$, that is, $D=1.45\lambda$, and compare the far field with BEAST EFIE solutions on meshes with $h=0.031$, which differ from those with $h=0.044$ by at most $0.04\%$. We mesh the surfaces quasi-uniformly, with $h$ from $0.25$ to $0.044$ and the equator made of mesh edges.

Fig.~\ref{fig:regularity} shows that $\varepsilon_F$ decreases with an approximate slope of $2$ on all surfaces down to $C^{1,0.5}$, even though the curvature of the $C^{1,0.5}$ surfaces is unbounded at the equator. On the Lipschitz surface with the concave edge, the error still converges, but more slowly, with a fitted slope of $1.7$ over all meshes. On the Lipschitz surface with the convex edge, it stagnates at $3.9\%$, and $\varepsilon_L$ stays at $0.035$ on all meshes. The computed field is thus not a Maxwell field. The iteration count is between $27$ and $37$ in all cases. Thus, reduced surface regularity does not affect the convergence down to $C^{1,0.5}$, and a concave edge, at which the field remains in $H^1$, only slows it. At a convex edge, where the field is not in $H^1$, the method fails, as anticipated in Section~I.

\subsection{Rounded Fichera Cube}\label{sec:fichera}
\begin{figure}[t]
\centering
\includegraphics{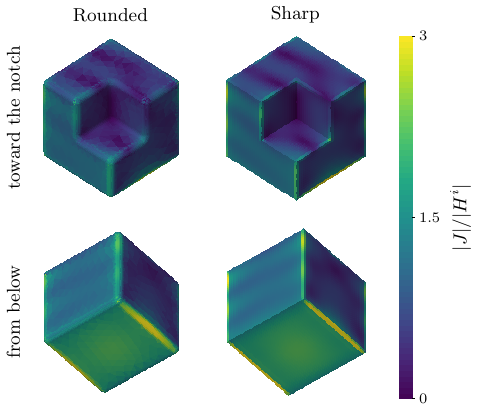}
\caption{Fichera cube with $D=2\lambda$, incidence along $+z$ with $E^i$ along $x$: $|J|/|H^i|$ on the rounded cube from the R-ECFOIE on the 7247-vertex mesh, recovered as in Section~\ref{sec:current_recovery}, and on the sharp cube from the BEAST EFIE on a uniform 6938-vertex mesh, on a common scale from $0$ to $3$, clipped above.}
\label{fig:fichera_current}
\end{figure}

The rounded Fichera cube, Fig.~\ref{fig:fichera_current}, is the cube $(-1,1)^3$ with the positive octant $[0,1)^3$ removed and all resulting edges rounded with fillet radius $\rho=0.1$, so that $D\approx3.32$. Since the fillets meet the faces tangentially, the surface is $C^{1,1}$: its normal is continuous, but its curvature jumps from $0$ on the faces to $1/\rho=10$ on the fillets. It thus falls outside the $C^{2,\alpha}$ setting of~\cite{burbano2026mahelmholtz}, but it has no edges and is smoother than the $C^{1,0.5}$ surfaces of Section~\ref{sec:regularity}, on which the errors still decrease at about second order. The Gmsh meshes, with $974$ to $73994$ vertices, are refined near the fillets according to the curvature. They are generated independently, so they are not nested, and their refinement near the fillets varies from mesh to mesh.

We first take $D=2\lambda$, that is, $k=\eta=4\pi/D\approx3.79$, and solve both equations without the constraints. The reference is a dense BEAST EFIE solution on a separately generated 2216-vertex mesh, which agrees to within $0.5\%$ in the far field with BEAST solutions on two other meshes.

\begin{table}[h]
\caption{Rounded Fichera Cube with $D=2\lambda$: Near- and Far-Field Errors and Longitudinal Fraction (\%) and GMRES Iterations Against the BEAST EFIE}
\label{tab:fichera_first}
\centering
\footnotesize
\setlength{\tabcolsep}{2.5pt}
\begin{tabular}{rrrrrrrrr}
\toprule
 & \multicolumn{4}{c}{R-ECFOIE} & \multicolumn{4}{c}{ECFOIE} \\
\cmidrule(lr){2-5}\cmidrule(lr){6-9}
$N$ & $\varepsilon_E$ & $\varepsilon_F$ & $\varepsilon_L$ & It. & $\varepsilon_E$ & $\varepsilon_F$ & $\varepsilon_L$ & It. \\
\midrule
974  & 17.64 & 14.91 & 9.07 & 59 & 13.05 & 11.05 & 6.55 & 215 \\
1569 & 12.96 & 11.08 & 6.33 & 63 & 9.92  & 8.42  & 4.90 & 250 \\
2486 & 8.54  & 7.35  & 4.01 & 67 & 5.62  & 4.72  & 2.85 & 288 \\
3206 & 8.74  & 7.42  & 4.28 & 68 & 6.08  & 5.07  & 3.15 & 308 \\
5094 & 5.50  & 4.68  & 2.64 & 71 & 3.70  & 3.09  & 1.89 & 388 \\
7247 & 3.83  & 3.26  & 1.84 & 73 & 2.58  & 2.16  & 1.32 & 427 \\
12975 & 2.12 & 1.78 & 1.06 & 73 & 1.39 & 1.15 & 0.74 & 512 \\
25887 & 1.14 & 0.95 & 0.59 & 74 & 0.83 & 0.69 & 0.44 & 623 \\
48089 & 0.64 & 0.54 & 0.31 & 73 & 0.54 & 0.48 & 0.23 & 740 \\
\bottomrule
\end{tabular}
\end{table}

Table~\ref{tab:fichera_first} shows that all errors and the longitudinal fraction decrease under refinement at about second order on the finer meshes, until the errors reach the level of the reference, and that the R-ECFOIE iteration count stays nearly constant. The ECFOIE is somewhat more accurate but needs four to ten times more iterations.

On the finest mesh, the recovered surface current differs from that of the BEAST EFIE on the same mesh by $2.5\%$ in the relative $L^2$ norm of $|J|$ on the flat faces and by $4.9\%$ overall, with the largest differences on the fillets. Fig.~\ref{fig:fichera_current} compares it with the EFIE current on the sharp cube. The two agree away from the edges, and on the lit face $z=-1$ both are close to the physical-optics value $2|H^i|$. At the edges, the current on the sharp cube concentrates in narrow lines, where it is singular, whereas the fillets spread it over a band.

\begin{table}[t]
\caption{Rounded Fichera Cube with $D=10^{-8}\lambda$: R-ECFOIE Near-Field Error $\varepsilon_E$ (\%) at Radii $d$ and $4d$, GMRES Iterations, and Induced Charge $q$, Without and With the Constraints}
\label{tab:fichera_low}
\centering
\footnotesize
\setlength{\tabcolsep}{2.5pt}
\begin{tabular}{rrrrcrrrc}
\toprule
 & \multicolumn{4}{c}{Without constraints} & \multicolumn{4}{c}{With constraints} \\
\cmidrule(lr){2-5}\cmidrule(lr){6-9}
$N$ & $d$ & $4d$ & It. & $q$ & $d$ & $4d$ & It. & $q$ \\
\midrule
974  & 80.80 & 291.5 & 77 & $3{\cdot}10^{-1}$ & 29.81 & 38.14 & 42 & $8{\cdot}10^{-6}$ \\
1569 & 21.76 & 48.39 & 80 & $2{\cdot}10^{-2}$ & 22.09 & 26.76 & 44 & $5{\cdot}10^{-6}$ \\
2486 & 27.84 & 105.8 & 83 & $1{\cdot}10^{-1}$ & 15.36 & 18.53 & 46 & $1{\cdot}10^{-5}$ \\
3206 & 14.37 & 26.72 & 82 & $2{\cdot}10^{-2}$ & 15.84 & 19.80 & 46 & $2{\cdot}10^{-5}$ \\
5094 & 8.639 & 20.25 & 84 & $5{\cdot}10^{-3}$ & 10.28 & 14.03 & 48 & $6{\cdot}10^{-6}$ \\
7247 & 9.272 & 32.84 & 87 & $3{\cdot}10^{-2}$ & 7.495 & 10.86 & 48 & $3{\cdot}10^{-6}$ \\
12975 & --- & --- & --- & --- & 4.382 & 6.875 & 49 & $4{\cdot}10^{-6}$ \\
25887 & --- & --- & --- & --- & 2.493 & 3.606 & 49 & $3{\cdot}10^{-6}$ \\
48089 & --- & --- & --- & --- & 1.337 & 1.879 & 50 & $2{\cdot}10^{-6}$ \\
73994 & --- & --- & --- & --- & 0.910 & 1.322 & 50 & $2{\cdot}10^{-6}$ \\
\bottomrule
\end{tabular}
\end{table}

At $D=10^{-8}\lambda$ we use incidence along $-x$ with $E^i$ along $y$, and on each mesh the reference is the electrostatic field computed on that mesh. Table~\ref{tab:fichera_low} shows the same behavior as for the sphere and the resonator, without and with the constraints, and with them GMRES needs about half as many iterations. On four finer meshes, solved only with the constraints, the errors continue to decrease at about second order in the mesh size, to below $1\%$ at radius $d$, with the same iteration count. The errors are nonetheless much larger than for the sphere and the Helmholtz resonator at comparable $N$, which shows that in both regimes the accuracy is limited by the resolution of the fillets. These results hold for the fillet radius $\rho=0.1$, and they do not extend to the sharp cube, whose edges the method does not handle. %

\subsection{Smoothed Stanford Bunny}\label{sec:bunny}
\begin{figure}[t]
\centering
\includegraphics{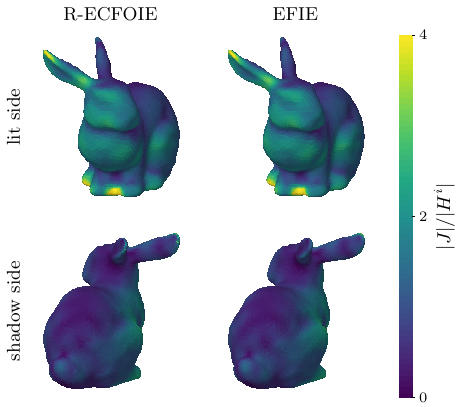}\\[2pt]
\includegraphics[width=\columnwidth]{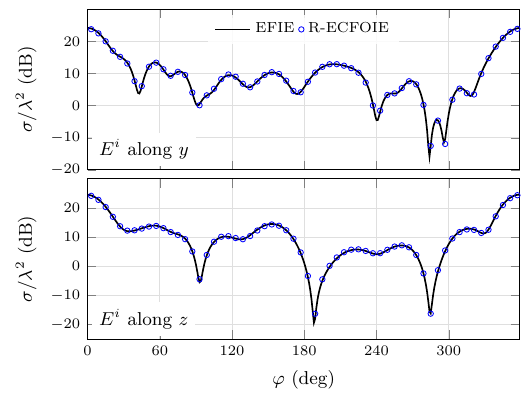}
\caption{Smoothed bunny at $k=\eta=4\pi$, incidence along $+x$, R-ECFOIE and BEAST EFIE on the same 9210-vertex mesh. Top: $|J|/|H^i|$ for $E^i$ along $y$, recovered as in Section~\ref{sec:current_recovery} for the R-ECFOIE, seen from the lit side ($-x$) and the shadow side, on a common scale from $0$ to $4$, clipped above. Bottom: bistatic RCS in the horizontal plane. Forward scattering is at $\varphi=0$.}
\label{fig:bunny_k4pi}
\end{figure}

Our last example is the Stanford bunny, a scanned geometry with strongly curved features. The scanned surface is not smooth: its triangulation has sharp folds at the ears and the base, which persist under refinement  and act like the convex edge of Section~\ref{sec:regularity}, where the R-ECFOIE fails. We therefore use a smoothed bunny, Fig.~\ref{fig:bunny_k4pi}, of diameter $D=2$, defined as a level set of the indicator function of the scan smoothed by Gaussian convolution. We remesh it with Gmsh and project the vertices onto the level set by Newton's method, so that refinement approximates the smooth surface. We then solve the problem of scattering of a plane wave propagating along $+x$, with $E^i$ along $y$ or along $z$.

At $k=\eta=2\pi$, against a BEAST EFIE reference on a separate 9104-vertex mesh, $\varepsilon_F$ decreases from $24.7\%$ and $22.0\%$ for $E^i$ along $y$ and $z$, respectively, on a 2289-vertex mesh to $7.1\%$ and $6.1\%$ on a 15624-vertex mesh.
At $k=\eta=4\pi$, where the bunny is four wavelengths across, we solve on a mesh with $9210$ vertices and $h\le\lambda/10=0.05$, and compare with the BEAST EFIE on the same mesh in Fig.~\ref{fig:bunny_k4pi}. The two surface currents show the same standing-wave pattern, and their magnitudes $|J|$ differ by $4.7\%$ in the relative $L^2$ norm, most at the tip of the upright ear. The far fields differ by $\varepsilon_F=4.0\%$ and $3.7\%$ for the two polarizations, and the R-ECFOIE needs $129$ GMRES iterations.

\section{Discussion and Conclusion}
We developed a continuous P1 Galerkin BEM for the electric field-only PEC formulation on flat triangular meshes, with a vertex-averaged curvature, Calder\'on regularization, and Lagrange-multiplier constraints for low frequencies. The errors decrease at about second order, also at $D=10^{-8}\lambda$, and the R-ECFOIE iteration count does not grow under refinement. On the rounded Fichera cube and the smoothed bunny, the far fields and surface currents agree with EFIE solutions.

The method is intended for smooth or rounded scatterers. This restriction comes from the formulation, which requires fields in $H^1$, and not from the discretization: the errors still decrease at about second order on $C^{1,\alpha}$ surfaces, but they stagnate on Lipschitz surfaces with convex edges. Open problems include a convergence analysis of the scheme, with a consistency estimate for the averaged curvature, and preconditioning near resonances.

\section*{Acknowledgment}
Claude (Anthropic) helped translate the first author's BEM code from Fortran~90 to Julia, implement geometry-specific features, and polish the text.

\bibliographystyle{IEEEtran}
\bibliography{references}

@article{burbano2026mahelmholtz,
  author  = {Burbano-Gallegos, Juan and P{\'e}rez-Arancibia, Carlos and Turc, Catalin},
  title   = {Maxwell {\`a} la {H}elmholtz: Electromagnetic Scattering by {3D} Perfect Electric Conductors via {H}elmholtz Integral Operators},
  journal = {ESAIM: Mathematical Modelling and Numerical Analysis},
  volume  = {60},
  number  = {1},
  pages   = {273--315},
  year    = {2026},
  doi     = {10.1051/m2an/2025101}
}

@article{rao1982scattering,
  author  = {Rao, Sadasiva M. and Wilton, Donald R. and Glisson, Allen W.},
  title   = {Electromagnetic Scattering by Surfaces of Arbitrary Shape},
  journal = {IEEE Transactions on Antennas and Propagation},
  volume  = {30},
  number  = {3},
  pages   = {409--418},
  year    = {1982},
  doi     = {10.1109/TAP.1982.1142818}
}

@article{costabeldauge2000singularities,
  author  = {Costabel, Martin and Dauge, Monique},
  title   = {Singularities of Electromagnetic Fields in Polyhedral Domains},
  journal = {Archive for Rational Mechanics and Analysis},
  volume  = {151},
  number  = {3},
  pages   = {221--276},
  year    = {2000},
  doi     = {10.1007/s002050050197}
}

@incollection{sullivan2008curvatures,
  author    = {Sullivan, John M.},
  title     = {Curvatures of Smooth and Discrete Surfaces},
  booktitle = {Discrete Differential Geometry},
  series    = {Oberwolfach Seminars},
  volume    = {38},
  publisher = {Birkh\"auser},
  year      = {2008},
  note      = {See also arXiv:0710.4497}
}

@book{sauter2011bem,
  author    = {Sauter, Stefan A. and Schwab, Christoph},
  title     = {Boundary Element Methods},
  series    = {Springer Series in Computational Mathematics},
  volume    = {39},
  publisher = {Springer},
  address   = {Berlin, Germany},
  year      = {2011},
  doi       = {10.1007/978-3-540-68093-2}
}

@article{geuzaine2009gmsh,
  author  = {Geuzaine, Christophe and Remacle, Jean-Fran{\c c}ois},
  title   = {{Gmsh}: A Three-Dimensional Finite Element Mesh Generator with Built-in Pre- and Post-Processing Facilities},
  journal = {International Journal for Numerical Methods in Engineering},
  volume  = {79},
  number  = {11},
  pages   = {1309--1331},
  year    = {2009},
  doi     = {10.1002/nme.2579}
}

@misc{hmatricesjl,
  author       = {{IntegralEquations contributors}},
  title        = {{HMatrices.jl}: Hierarchical Matrices for Julia},
  howpublished = {\url{https://github.com/IntegralEquations/HMatrices.jl}},
  year         = {2026}
}

@misc{beastjl,
  author       = {{BEAST.jl contributors}},
  title        = {{BEAST.jl}: Boundary Element Analysis and Simulation Toolkit},
  howpublished = {\url{https://github.com/krcools/BEAST.jl}},
  year         = {2026}
}

@article{hildebrandt2006convergence,
  author  = {Hildebrandt, Klaus and Polthier, Konrad and Wardetzky, Max},
  title   = {On the Convergence of Metric and Geometric Properties of Polyhedral Surfaces},
  journal = {Geometriae Dedicata},
  volume  = {123},
  pages   = {89--112},
  year    = {2006},
  doi     = {10.1007/s10711-006-9109-5}
}

@book{steinbach2008numerical,
  author    = {Steinbach, Olaf},
  title     = {Numerical Approximation Methods for Elliptic Boundary Value Problems: Finite and Boundary Elements},
  publisher = {Springer},
  address   = {New York},
  year      = {2008},
  doi       = {10.1007/978-0-387-68805-3}
}

@incollection{meyer2003discrete,
  author    = {Meyer, Mark and Desbrun, Mathieu and Schr{\"o}der, Peter and Barr, Alan H.},
  title     = {Discrete Differential-Geometry Operators for Triangulated 2-Manifolds},
  booktitle = {Visualization and Mathematics III},
  editor    = {Hege, Hans-Christian and Polthier, Konrad},
  publisher = {Springer},
  address   = {Berlin},
  pages     = {35--57},
  year      = {2003},
  doi       = {10.1007/978-3-662-05105-4_2}
}

@article{wilton1984potential,
  author  = {Wilton, D. R. and Rao, S. M. and Glisson, A. W. and Schaubert, D. H. and Al-Bundak, O. M. and Butler, C. M.},
  title   = {Potential Integrals for Uniform and Linear Source Distributions on Polygonal and Polyhedral Domains},
  journal = {IEEE Transactions on Antennas and Propagation},
  volume  = {32},
  number  = {3},
  pages   = {276--281},
  year    = {1984},
  doi     = {10.1109/TAP.1984.1143304}
}

@book{coltonkress1983integral,
  author    = {Colton, David and Kress, Rainer},
  title     = {Integral Equation Methods in Scattering Theory},
  publisher = {SIAM},
  address   = {Philadelphia, PA},
  year      = {2013},
  note      = {Reprint of the 1983 Wiley edition},
  doi       = {10.1137/1.9781611973167}
}

@book{volakis2012integral,
  author    = {Volakis, John L. and Sertel, Kubilay},
  title     = {Integral Equation Methods for Electromagnetics},
  publisher = {SciTech Publishing},
  year      = {2012},
  doi       = {10.1049/SBEW045E}
}

@article{bebendorf2000approximation,
  author  = {Bebendorf, Mario},
  title   = {Approximation of Boundary Element Matrices},
  journal = {Numerische Mathematik},
  volume  = {86},
  number  = {4},
  pages   = {565--589},
  year    = {2000},
  doi     = {10.1007/s002110000192}
}

@incollection{buffa2003galerkin,
  author    = {Buffa, Annalisa and Hiptmair, Ralf},
  title     = {Galerkin Boundary Element Methods for Electromagnetic Scattering},
  booktitle = {Topics in Computational Wave Propagation: Direct and Inverse Problems},
  editor    = {Ainsworth, Mark and Davies, Penny and Duncan, Duncan and Martin, Paul and Rynne, Brendan},
  series    = {Lecture Notes in Computational Science and Engineering},
  volume    = {31},
  pages     = {83--124},
  publisher = {Springer},
  address   = {Berlin, Germany},
  year      = {2003},
  doi       = {10.1007/978-3-642-55483-4_3}
}

@article{buffa2007dual,
  author  = {Buffa, Annalisa and Christiansen, Snorre H.},
  title   = {A Dual Finite Element Complex on the Barycentric Refinement},
  journal = {Mathematics of Computation},
  volume  = {76},
  number  = {260},
  pages   = {1743--1769},
  year    = {2007},
  doi     = {10.1090/S0025-5718-07-01965-5}
}

@article{zhao2000integral,
  author  = {Zhao, Jun-Sheng and Chew, Weng Cho},
  title   = {Integral Equation Solution of {Maxwell}'s Equations from Zero Frequency to Microwave Frequencies},
  journal = {IEEE Transactions on Antennas and Propagation},
  volume  = {48},
  number  = {10},
  pages   = {1635--1645},
  year    = {2000},
  doi     = {10.1109/8.899680}
}

@incollection{merlini2020lowfrequency,
  author    = {Merlini, Adrien and D{\'e}ly, Alexandre and Cools, Kristof and Andriulli, Francesco P.},
  title     = {Electromagnetic Modelling at Arbitrarily Low Frequency via the Quasi-{Helmholtz} Projectors},
  booktitle = {Advances in Mathematical Methods for Electromagnetics},
  editor    = {Kobayashi, Kazuya and Smith, Paul D.},
  pages     = {381--415},
  publisher = {Institution of Engineering and Technology},
  year      = {2020},
  doi       = {10.1049/SBEW528E_ch16}
}

@article{andriulli2008multiplicative,
  author  = {Andriulli, Francesco P. and Cools, Kristof and Ba{\u{g}}c{\i}, Hakan and Olyslager, Femke and Buffa, Annalisa and Christiansen, Snorre and Michielssen, Eric},
  title   = {A Multiplicative {C}alder{\'o}n Preconditioner for the Electric Field Integral Equation},
  journal = {IEEE Transactions on Antennas and Propagation},
  volume  = {56},
  pages   = {2398--2412},
  year    = {2008}
}

@article{taskinen2006current,
  author  = {Taskinen, Matti and Yl{\"a}-Oijala, Pasi},
  title   = {Current and Charge Integral Equation Formulation},
  journal = {IEEE Transactions on Antennas and Propagation},
  volume  = {54},
  pages   = {58--67},
  year    = {2006}
}

@article{epstein2010debye,
  author  = {Epstein, Charles L. and Greengard, Leslie},
  title   = {Debye Sources and the Numerical Solution of the Time Harmonic {M}axwell Equations},
  journal = {Communications on Pure and Applied Mathematics},
  volume  = {63},
  pages   = {413--463},
  year    = {2010}
}

@article{ganesh2006spectrally,
  author  = {Ganesh, M. and Hawkins, S. C.},
  title   = {A Spectrally Accurate Algorithm for Electromagnetic Scattering in Three Dimensions},
  journal = {Numerical Algorithms},
  volume  = {43},
  pages   = {25--60},
  year    = {2006}
}

@article{bruno2009electromagnetic,
  author  = {Bruno, Oscar and Elling, Tim and Paffenroth, Randy and Turc, Catalin},
  title   = {Electromagnetic Integral Equations Requiring Small Numbers of {K}rylov-Subspace Iterations},
  journal = {Journal of Computational Physics},
  volume  = {228},
  pages   = {6169--6183},
  year    = {2009}
}

@article{faria2021general,
  author  = {Faria, Luiz M. and P{\'e}rez-Arancibia, Carlos and Bonnet, Marc},
  title   = {General-Purpose Kernel Regularization of Boundary Integral Equations via Density Interpolation},
  journal = {Computer Methods in Applied Mechanics and Engineering},
  volume  = {378},
  pages   = {113703},
  year    = {2021},
  doi     = {10.1016/j.cma.2021.113703}
}

@article{hu2021chebyshev,
  author  = {Hu, Jin and Garza, Emmanuel and Sideris, Constantine},
  title   = {A {C}hebyshev-Based High-Order-Accurate Integral Equation Solver for {M}axwell's Equations},
  journal = {IEEE Transactions on Antennas and Propagation},
  volume  = {69},
  pages   = {5790--5800},
  year    = {2021}
}

@article{hu2024continuity,
  author  = {Hu, Jin and Sideris, Constantine},
  title   = {A High-Order {Nystr{\"o}m}-Based Scheme Explicitly Enforcing Surface Density Continuity for the Electric Field Integral Equation},
  journal = {IEEE Antennas and Wireless Propagation Letters},
  volume  = {23},
  number  = {9},
  pages   = {2633--2637},
  year    = {2024},
  doi     = {10.1109/LAWP.2024.3402749}
}

@misc{hofmann2026continuity,
  author        = {Hofmann, Bernd and Molavi, Reza and Sideris, Constantine},
  title         = {High-Order-Accurate Continuity Enforcing {Nystr{\"o}m} Discretization of {3D} {Maxwell} Combined Field Integral Equations},
  year          = {2026},
  eprint        = {2608.28876},
  archivePrefix = {arXiv},
  primaryClass  = {math.NA}
}

@article{yuffa2018tensorial,
  author  = {Yuffa, Alex J. and Markkanen, Johannes},
  title   = {A {3-D} Tensorial Integral Formulation of Scattering Containing Intriguing Relations},
  journal = {IEEE Transactions on Antennas and Propagation},
  volume  = {66},
  pages   = {5274--5281},
  year    = {2018}
}

@article{klaseboer2017nonsingular,
  author  = {Klaseboer, Evert and Sun, Qiang and Chan, Derek Y. C.},
  title   = {Nonsingular Field-Only Surface Integral Equations for Electromagnetic Scattering},
  journal = {IEEE Transactions on Antennas and Propagation},
  volume  = {65},
  number  = {2},
  pages   = {972--977},
  year    = {2017}
}

@article{sun2020pec,
  author  = {Sun, Qiang and Klaseboer, Evert and Yuffa, Alex J. and Chan, Derek Y. C.},
  title   = {Field-Only Surface Integral Equations: Scattering from a Perfect Electric Conductor},
  journal = {Journal of the Optical Society of America A},
  volume  = {37},
  pages   = {276--283},
  year    = {2020}
}

@article{faria2024combined,
  author  = {Faria, Luiz M. and P{\'e}rez-Arancibia, Carlos and Turc, Catalin},
  title   = {Combined Field-Only Boundary Integral Equations for {PEC} Electromagnetic Scattering Problem in Spherical Geometries},
  journal = {SIAM Journal on Applied Mathematics},
  volume  = {84},
  pages   = {19--38},
  year    = {2024}
}

@misc{perez2026direct,
  author       = {P{\'e}rez-Arancibia, Carlos and Turc, Catalin},
  title        = {Maxwell {\`a} la {H}elmholtz: Direct Boundary Integral Equations for {3D} Scattering by Perfect Electric Conductors via {H}elmholtz Operators},
  howpublished = {arXiv:2605.01670},
  year         = {2026}
}

@article{bonnetbendhia1999singular,
  author  = {Bonnet-Ben Dhia, Anne-Sophie and Hazard, Christophe and Lohrengel, St{\'e}phanie},
  title   = {A Singular Field Method for the Solution of {M}axwell's Equations in Polyhedral Domains},
  journal = {SIAM Journal on Applied Mathematics},
  volume  = {59},
  number  = {6},
  pages   = {2028--2044},
  year    = {1999},
  doi     = {10.1137/S0036139997323383}
}

@book{dassios2000low,
  author    = {Dassios, George and Kleinman, Ralph},
  title     = {Low Frequency Scattering},
  publisher = {Clarendon Press},
  address   = {Oxford},
  series    = {Oxford Mathematical Monographs},
  year      = {2000}
}

@article{marchand2022applying,
  author  = {Marchand, Pierre and Galkowski, Jeffrey and Spence, Alastair and Spence, Euan A.},
  title   = {Applying {GMRES} to the {H}elmholtz equation with strong trapping: how does the number of iterations depend on the frequency?},
  journal = {Advances in Computational Mathematics},
  volume  = {48},
  number  = {4},
  pages   = {37},
  year    = {2022},
  doi     = {10.1007/s10444-022-09931-9}
}

@article{boubendir2013wave,
  author  = {Boubendir, Yassine and Turc, Catalin},
  title   = {Wave-number estimates for regularized combined field boundary integral operators in acoustic scattering problems with {N}eumann boundary conditions},
  journal = {IMA J. Numer. Anal.},
  volume  = {33},
  number  = {4},
  pages   = {1176--1225},
  year    = {2013}
}

\end{document}